\documentclass[a4paper,10pt]{amsart}
\usepackage[all]{xy}
\usepackage{amsfonts}
\usepackage[centertags]{amsmath}
\usepackage{amssymb}
\usepackage{amsthm}
\usepackage{amscd}
\newtheorem{statement}{\bf Statement}

\newtheorem{corollary}{Corollary}[statement]

\theoremstyle{remark}

\newtheorem{example}{Example}

\title{Properties of two selections in metric spaces of Busemann  nonpositive curvature}
\author{Pyotr N. Ivanshin}

\begin{document}

\maketitle

\begin{abstract}
In this article the author presents results on the selection for the space of convex hulls of $n$ points and compacts of the complete convex metric space of Busemann nonpositve curvature. Namely, we determine Lipschitz and H\"older properties of barycenter and cirumcenter mappings.
\end{abstract}

\section{Introduction}
In this article the author presents results on the selection for the space of convex hulls of $n$ points and compacts of the complete metric space. This problem seems to be of interest for some time and has the abundant set of solutions in Hadamard spaces \cite{LPS, SB}, in Banach spaces \cite{Shvarts} and in special metric spaces \cite{Khamsi1}. The same can be said about properties of weighted points of the finite sets \cite{Vesely, LPS}. Note also that there exists the description of the Chebyshev center behaviour \cite{SV, BI} in special Banach spaces. Note that here we present the result generalising Theorem 1 of the latter paper to nonlinear spaces. Moreover we show that the behaviour of the distance function $d(\mathrm{cheb}(A), \mathrm{cheb}(B))$, especially non-Lipschitz side of it, is controlled only by convexity properties of the unit ball of metric space. 

\subsection{Notation} Here we mostly study sets of two spaces, namely $\Sigma_n(H)$, being the space of all $n$-nets of $H$ and $K(H)$ --- space of all compact convex subsets of $H$.

Also we denote by $\mathrm{cheb}(M)$ the set of Chebyshev centers of $M \subset H$ and by $\mathrm{diam}(M)$ diameter of the set $M$. Also $\mathrm{Hd}(V, W)$ equals Hausdorff distance between $V, W \subset H$.

Note also that we usually assume that space $H$ is such that any two points can be connected with unique geodesic line.

\subsection{Main results} In the paper we present inductive construction of the desired point. There are two approaches using 1) mean point (in some sense generalising barycenter) or 2) Chebyshev center of the set. At first we investigate properties of these mappings in non-linear spaces. The first of them is correctness of the definition of the mean point.

{\bf Statement 1.}{\it
There exists a common limit point $\mathrm{mp}(\sigma)$ for all sequences of points $(x^k_i)_{k=1}^{\infty}$, $i =1, \ldots, n$.}

Then we must consider continuity properties of this mapping.

Denote by $x_1, x_2, \ldots , x_n \in H$ elements of $\sigma \in \Sigma_n(N)$. Then the following holds true.

{\bf Statement 2.}{\it 
The mapping $\mathrm{mp}: (\Sigma_n(H), \mathrm{Hd}) \to H$ is $1$-Lipshitz in $\min\{d(x_i, x_j)| i, j =1, \ldots, n, i\neq j \}/2$-neighbourhood of $\sigma$ .}

Moreover the following holds true:

{\bf Statement 7.}{\it 
The mean point exists for any $V \in K(H)$.}

We note also that $\mp(V)$ coincides with usual barycenter of $V$ in case $H$ is Hadamard space. 

Then we apply methods of \cite{IS, IS2} to prove the following statement:

{\bf Statement 9.}{\it 
The Chebyshev selection map $\mathrm{cheb}: K(H) \to H$, $V \mapsto \mathrm{cheb}(V)$ is generalised H\"older with power constant $1/2$ \cite{BI, SV} in the Hadamard space $H$.}

{\bf Statement 10.}{\it
Let $H$ be strictly convex space of Busemann nonpositive curvature.
Then the upper bound for the power of the generalised H\"older map $\mathrm{cheb}: B(H) \to H$ is $1/2$.}

Finally, we construct Lipschitz selection points for the convex hulls of no more than $n$ points of $H$. Naturally, Lipschitz constants depends on the number of points, since, for instance, there can be no Lipschitz selection in this case even in Banach spaces \cite{PY, Lang}.

\section{Mean point.}

\subsection{Mean point of $n$-net.}
Let $H$ be a Busemann nonpositively curved space.

Let $\sigma \in \Sigma_n(H)$, $\sigma=\{x_1, \ldots, x_n\}$ such that $x_i\neq x_j$ for $i \neq j$.
Consider the inductive construction of the mean point of $\sigma$. For $n=2$ put $\mathrm{mp}(\{x_1, x_2\})=m(x_1, x_2)$. Assume now that we have correctly defined mapping $\mathrm{mp}: \Sigma_{n-1}(H) \to H$. Then for $\sigma=\{x_1, \ldots, x_n\}$ one may consider the set $\sigma^1=\{x^1_1, \ldots, x_n^1\}$, here $x_i^1= \mathrm{mp}\{x_1, \ldots, \widehat{x_i}, \ldots, x_n\}$ and $\widehat{x_i}$ means that this point is excluded from $\sigma$. Then repeat the procedure with $\sigma^1$ to get $\sigma^2$.

\begin{statement}
There exists a common limit point $\mathrm{mp}(\sigma)$ for all sequences of points $(x^k_i)_{k=1}^{\infty}$, $i =1, \ldots, n$. 
\end{statement}

\begin{proof}
It suffices to prove that $d(x^k_i, x^k_j) \leq 1/2 d(x^{k-1}_i, x^{k-1}_j)$. 

The proof is by induction on $n$ (the number of points). 
The base of the induction is the first consistent case $n=3$. The statement follows from the definition of the nonpositive curvature. 
Now assume that the statement holds true for all $l \leq n-1$ and prove it for $n$. Let us prove it for $\sigma$ and $\sigma^1$ because the rest will easily follow from it. Fix any pair of points $x_i, x_j \in \sigma$.
Recall the construction of $x^1_i, x^1_j$ and note that we start from the mean points of the sets $\{x_i, x_1, \ldots,\widehat{x_i}, \ldots, \widehat{x_i}, \ldots, \widehat{x_l}, x_n\}$ and $\{x_j, x_1, \ldots,\widehat{x_i}, \ldots, \widehat{x_j}, \ldots, \widehat{x_k}, x_n\}$. For each pair of points $\mathrm{mp}(\{x_i, x_1, \ldots,\widehat{x_i}, \ldots, \widehat{x_i}, \ldots, \widehat{x_l}, x_n\})$ and $\mathrm{mp}(\{x_j, x_1, \ldots,\widehat{x_i}, \ldots, \widehat{x_j}, \ldots, \widehat{x_k}, x_n\})$ from the described sets we have the desired inequality for the triple $x_i, x_j, x_k$ and since by induction hypothesis points of $\{x_i, x_1, \ldots,\widehat{x_i}, \ldots, \widehat{x_i}, \ldots, \widehat{x_l}, x_n\}$ and $\{x_j, x_1, \ldots,\widehat{x_i}, \ldots, \widehat{x_j}, \ldots, \widehat{x_k}, x_n\}$ converge to the limit points $x^1_i, x^1_j$ respectively, we get the inequality.
\end{proof}

\begin{statement}
The mapping $\mathrm{mp}: (\Sigma_n(H), \mathrm{Hd}) \to H$ is $1$-Lipshitz in $\min\limits_{x_i \neq x_j \in \sigma}{d(x_i, x_j)}/2$-neighbourhood of each $\sigma \in \Sigma_n(H)$.
\end{statement}

\begin{proof}
The proof is again by induction.
The base of the induction is $\Sigma_2(H)$ for which the inequality follows from the convexity of the distance function between two segments.

The rest of the proof repeats corresponding part of the proof of Statement 1.
\end{proof}

Let $\sigma \in \Sigma_n(H)$, $\sigma=\{x_1, \ldots, x_n\}$.
Recall now the definition of the barycenter $\mathrm{bc}(\sigma)$ as the point providing minimum to the function $f_{\sigma}(x)=\sum\limits_{i=1}^{n} d^2(x, x_i)$.
Since in Busemann space the relation $d^2(z, \gamma(t)) \leq (1-t) d^2(z, \gamma(0))+ t d^2(z, \gamma(1))- t(1-t)d^2(\gamma(0), \gamma(1))$ does not necessarily take place, we can not apply arguments of \cite{LPS}; nevertheless we get the following proposition:

\begin{statement}
The mean point $\mathrm{mp}(\sigma)$ minimizes function which is majorized by barycenter one for any $\sigma \in \Sigma_n(H)$ in Hadamard space $H$.
\end{statement}

\begin{proof}
The proof is again by induction on $n$. The base of the induction is the trivial case of $n=2$. Then $\mathrm{bc}(\{x_1, x_2\})=m(x_1, x_2)=\mathrm{mp}(\{x_1, x_2\})$.

It suffices to prove that $\forall x \in \mathrm{co}(\sigma) \setminus \sigma^1$, $f_{\sigma}(x) \geq \sum\limits_{i=1}^{n} d^2(x_i, \mathrm{co}(\sigma^1))$ since then the first approximation of the point $\mathrm{mp}(\sigma)$ is also the first approximation of the barycenter.

In order to prove it one must consider the auxiliary number $\frac{1}{n-1} \sum\limits_{i=1}^{n}\sum\limits_{j\neq i} d^2(x_i, x^1_j)$. Then convexity of both $\mathrm{co}(\sigma^1)$ and distance function implies that $\sum\limits_{i=1}^{n} d^2(x_i, \mathrm{co}(\sigma^1)) \leq \frac{1}{n-1} \sum\limits_{i=1}^{n}\sum\limits_{j\neq i} d^2(x_i, x^1_j)$. At the same time for any $x \in H$ we have the relation $\sum\limits_{i=1}^{n} d^2(x, x_i) \geq \frac{1}{n-1} \sum\limits_{i=1}^{n}\sum\limits_{j\neq i} d^2(x_i, x^1_j)$. Combination of these two inequalities completes the proof for the set $\sigma^1$.

Since the space $H$ is Hadamard the "parallelogramm inequality" $(n-1)d^2(x, x^1_i) \leq \sum\limits_{j=1, j \neq i}^{n}d^2(x, x_j)-\sum\limits_{j=1, j \neq i}^{n}d^2(x^1_i, x_j)$ holds true. Note than in Euclidean space it becomes an equality. Thus $(n-1)\sum\limits_{i=1}^{n}d^2(x, x^1_i)+\sum\limits_{i=1}^{n}\sum\limits_{j=1, j \neq i}^{n}d^2(x^1_i, x_j) \leq \sum\limits_{i=1}^{n}\sum\limits_{j=1, j \neq i}^{n}d^2(x, x_j)=(n-1) \sum\limits_{j=1}^{n}d^2(x, x_j)$. Here the second summond $\sum\limits_{i=1}^{n}\sum\limits_{j=1, j \neq i}^{n}d^2(x^1_i, x_j)$ is minimised in the first step of the proof and the first summond $\sum\limits_{i=1}^{n}d^2(x, x^1_i)$ can be analysed similarly. This completes the proof.
\end{proof}

Note that there is no such relation between mean point and barycenter of the set in the space of Busemann nonpositive curvature.

\begin{example}
If $H =(\mathbb{R}^n, \|\cdot\|_p)$ or $L^p(\mathbb{R})$ for $p\geq2$ then for any $\sigma \in \Sigma_k(H)$ the point $\mathrm{mp}(\sigma)$ minimises the function $f: H \to \mathbb{R}^+$, $f(x)=\sum\limits_{j=1}^{k} \|x - x_j\|_p^p$.

Let us first consider the simplest nontrivial case of $(\mathrm{R}^n, \|\cdot\|_p)$ and $\Sigma_3(\mathrm{R}^n)$. Then the statement follows from the inequality $1+t^p \leq (1+t)^p + (1-t)^p$ for $t \in (0,1)$.

Now consider the case of $\Sigma_{k+1}(\mathrm{R}^n)$. Then for arbitrary point $x \in H$ 
$$\sum\limits_{i=1}^{k} \|x-x_i\|^p=\sum\limits_{i=1}^{k} \|pr(x)-x_i\|^p+k \|x-pr(x)\|^p= 
$$
$$
\sum\limits_{i=1}^{k} \|x^1_{k+1} + (pr(x)-x^1_{k+1})-x_i\|^p+k \|x-pr(x)\|^p \geq
$$
$$
\geq \sum\limits_{i=1}^{k} \|x^1_{k+1}-x_i\|^p + k \|x-pr(x)\|^p+
p \sum\limits_{i=1}^{k}\sum\limits_{j=1}^{n}|(x_i)_j|^{p-1} (pr(x)-x^1_{k+1})_j+ 
$$
$$
+\frac{p (p-1)}{2}\sum\limits_{i=1}^{k}\sum\limits_{j=1}^{n}|(x_i)_j|^{p-2} |(pr(x)-x^1_{k+1})_j|^2 + \ldots = 
$$
$$
=\sum\limits_{i=1}^{k} \|x^1_{k+1}-x_i\|^p + k \|x-pr(x)\|^p + k\| pr(x)-x^1_{k+1} \|^p  - k\| pr(x)-x^1_{k+1} \|^p + 
$$
$$
+p \sum\limits_{i=1}^{k}\sum\limits_{j=1}^{n}|(x_i)_j|^{p-1} (pr(x)-x^1_{k+1})_j+ \frac{p (p-1)}{2}\sum\limits_{i=1}^{k}\sum\limits_{j=1}^{n}|(x_i)_j|^{p-2} |(pr(x)-x^1_{k+1})_j|^2 + \ldots.
$$
Hence it suffices to prove that $- k\| pr(x)-x^1_{k+1} \|^p + p \sum\limits_{i=1}^{k}\sum\limits_{j=1}^{n}|(x_i)_j|^{p-1} (pr(x)-x^1_{k+1})_j+ \frac{p (p-1)}{2}\sum\limits_{i=1}^{k}\sum\limits_{j=1}^{n}|(x_i)_j|^{p-2} |(pr(x)-x^1_{k+1})_j|^2 \geq 0$. Since $x^1_{k+1}$ is the minimum point of $\sum\limits_{i=1}^{k}\|x -x_i\|^p$  $\sum\limits_{i=1}^{k}\sum\limits_{j=1}^{n}|(x_i)_j|^{p-1} (pr(x)-x^1_{k+1})_j=0$ and since $l \in \mathrm{co}\{x_1, \ldots, n\}$ $\sum\limits_{j=1}^{n}|(x_i)_j|^{p-2}>\|l_i\|^{p-2}$ for any $i=1, \ldots n$. This completes the proof.
\end{example}

\subsection{Mean point construction for the space of nonpositive Busemann curvature.}

Let us define a mean point of the arbitrary set $V \subset H$. Consider the sequence $\sigma_n \in \Sigma_n(H)$, $\sigma_n \subset V$ as in the previous section. Then call the limit point of the sequence (in case such a point exists) $\mathrm{mp} (\sigma_n)$ the mean point of $V$.

\begin{statement}
The mean point exists for any compact infinite set $V =\{x_1, x_2, \ldots| x_i \in H \}$ such that $\sum\limits_{i=1}^{\infty}d(x_i, x_{i+1})$ converges.
\end{statement}

\begin{proof}
Construction of the mean point implies that for any $\sigma \subset V$ $d(\mathrm{mp}(\sigma), \mathrm{mp}(\sigma \bigcup\limits_{j=1}^{k}\{x_j\})) \leq \max\{d(x_j, \mathrm{co}(\sigma))| j= 1, \ldots k\}$.

First let us show that for $\sigma \in \Sigma_n(H)$ $\sigma=\{x_1, \ldots, x_n\}$ $d(\mathrm{mp}(\sigma), \mathrm{mp}(\sigma \bigcup\{x\})) \leq \max\limits_{i=1, \ldots, n}\{d(x, x_i)\}$. The proof is by induction on $n$. The base of the induction is as usual $n=1$, in which case the proof is trivial. 
So assume that for $k \leq n-1$ the claim holds true and prove it for $n$. So we must estimate $d(\mathrm{mp}(\sigma), \mathrm{mp}(\sigma \bigcup\{x\}))$. Recall the construction of $\mathrm{mp}(\sigma)$. Then any $n$-subnet of $\sigma \bigcup\{x\}$ is either $\sigma$ or $\delta_{n-1} \bigcup \{x\}$, here $\delta_{n-1}$ is some $n-1$-subnet of $\sigma$. By induction hypothesis the distances $d(\mathrm{mp}(\delta_{n-1} \bigcup \{x\}), \mathrm{mp}(\delta_{n-1}))$ satisfy the desired relation. Hence the result follows from the convexity of the distance function between sets.

Now consider the general case of $\sigma$ and $\sigma\bigcup\limits_{j=1}^{k}\{x_j\}$. 
\end{proof}

\begin{statement}
Let $\sigma \in \Sigma_{n-1}$, $\sigma=\{x_1, \ldots, x_{n-1}\}$ and $x \in H$ then $d(\mathrm{mp}(\sigma_n), \mathrm{mp}(\sigma \bigcup \{x\})) \leq \max\limits_{i=1, \ldots, n-1}\{d(x, x_i)\}/n$.
\end{statement}

\begin{proof}
The proof is again by induction on $n$. Assume that the statement holds true for $\Sigma_{n-2}$. Then recall the construction from the proof of Statement ??. Thus $d(\mathrm{mp}(\sigma), \mathrm{mp}(\sigma \bigcup \{x\})) \leq (\frac{1}{n-1}-\frac{1}{n-1}(\frac{1}{n-1})+\frac{1}{n-1} (\frac{1}{n-1}(\frac{1}{n-1}))+ \ldots ) \max\limits_{i=1, \ldots, n}\{d(x, x_i)\}=(\frac{1}{n-1}-\frac{1}{(n-1)^2}+\frac{1}{(n-1)^3}-\ldots ) \max\limits_{i=1, \ldots, n}\{d(x, x_i)\}=\frac{\max\limits_{i=1, \ldots, n}\{d(x, x_i)\}}{n}$.
\end{proof}

Similar considerations enable us to get

\begin{statement}
Let $\sigma \in \Sigma_n$, $\sigma=\{x_1, \ldots, x_n\}$ and $\sigma' \in \Sigma_k$, $\sigma'=\{y_1, \ldots, y_k\}$ then $d(\mathrm{mp}(\sigma_n), \mathrm{mp}(\sigma_n \bigcup \sigma')) \leq \max\limits_{i=1, \ldots, n, j=1, \ldots, k}\{d(x_i, y_j)\} \frac{k}{n+k}$.
\end{statement}

\begin{proof}
The proof is by induction on $k+n$. As in the previous statement we get an estimate $d(\mathrm{mp}(\sigma), \mathrm{mp}(\sigma \bigcup \sigma')) \leq (\frac{k}{n+k-1}-\frac{1}{n+k-1}(\frac{k}{n+k-1})+\frac{1}{n+k-1} (\frac{1}{n+k-1}(\frac{k}{n+k-1}))+ \ldots ) \max\limits_{i=1, \ldots, n}\{d(x, x_i)\}=(\frac{k}{n+k-1}-\frac{k}{(n+k-1)^2}+\frac{k}{(n+k-1)^3}-\ldots ) \max\limits_{i=1, \ldots , n}\{d(x, x_i)\}=\max\limits_{i=1, \ldots, n, j=1, \ldots, k}\{d( x_i, y_j)\}\frac{k }{n+k}$.
\end{proof}

Assume now that we vary masses of the points so that $\sigma= \bigcup\limits_{i=1}^{n} \{x_i, m_i\}$ and $\sigma'=\bigcup\limits_{i=1}^{n} \{x_i, m^{'}_i\}$, here $m_i, m^{'}_i \geq 0$, $i= 1, \ldots , n$.

\begin{corollary}
$d(\mathrm{mp}(\sigma), \mathrm{mp}(\sigma')) \leq \mathrm{diam}(\sigma)\sum\limits_{i=1}^{n}|m_i -m^{'}_{i}|$.
\end{corollary}

\begin{proof}
Assume first that $m_i, m^{'}_i \in \mathbb{Q}^+$. Hence $m_i=\frac{p_i}{q_i}$ and $m^{'}_i=\frac{p^{'}_i}{q^{'}_i}$.  Consider $q=\mathrm{LCM}\{q_i, q^{'}_i | i=1, \ldots, n\}$.
Now we must represent $x_i \in \sigma$ and $x^{'}_{i} \in \sigma'$ as the union of $p_i q/q_i$ and $p^{'}_i q/q^{'}_i$ of points of the same mass $1/q$, respectively. 

Then we apply statement 6 to get the result.
The general case easily follows.
\end{proof}

\begin{statement}
The mean point exists for any $V \in K(H)$.
\end{statement}

\begin{proof}
Fix $\varepsilon>0$.
Consider an $\varepsilon$-net $\sigma_1$ of $V$ such that the measure of the set $\sigma_1^{'}=\{x \in V| \exists i, j, i \neq j, d(x, x_i) <\varepsilon,  d(x, x_j) <\varepsilon\} < \delta \mu(V)$. Then for any other uniformly distributed $\varepsilon$-net $\sigma_2$ the set of points $\sigma_2^{'}$ for which we  have no one-to-one correspondence between $\sigma_1$ and $\sigma_2$ consists of no more than $\delta |\sigma_2|$ points if cardinalities of $\sigma_1$ and $\sigma_2$ coincide.
Thus statement 6 and its corollary imply that $d(\mathrm{mp}(\sigma_1), \mathrm{mp}(\sigma_2)) \leq d(\mathrm{mp}(\sigma_1), \mathrm{mp}(\sigma_1 \setminus \sigma_1^{'}))+d(\mathrm{mp}(\sigma_1\setminus \sigma_1^{'}), \mathrm{mp}(\sigma_2\setminus \sigma_2^{'}))+ d(\mathrm{mp}(\sigma_2), \mathrm{mp}(\sigma_1 \setminus \sigma_2^{'})) < 2 \delta + 2 \varepsilon$. This completes the proof.
\end{proof}

Evidently we arrive to the statement analogous to one of \cite{LPS}.

\begin{corollary}
Let $\mu_1$ and $\mu_2$ be two probability measures with supports $V_1$ and $V_2$, absolutely continuous with respect to the Hausdorff measure $\mu_H$ with Radon-Nikodym derivatives $\theta_1$ and $\theta_2$. Then $d(\mathrm{mp}(V_1), \mathrm{mp}(V_2)) \leq \int\limits_{V_1 \bigcup V_2} |\theta_1 -\theta_2| d\mu_H$.
\end{corollary}

\begin{proof}
It suffices to prove the statement for the finite sets $\sigma_1$ and $\sigma_2$ consisting of points of equal mass. Consider $\sigma=\sigma_1 \bigcap \sigma_2$. Then triangle inequality combined with statement 7 implies that
$d(\mathrm{mp}(\sigma_1), \mathrm{mp}(\sigma_2)) \leq d(\mathrm{mp}(\sigma_1), \mathrm{mp}(\sigma))+ d(\mathrm{mp}(\sigma_2), \mathrm{mp}(\sigma)) \leq \mathrm{diam}(\sigma_1 \bigcup \sigma_2) (\frac{|\sigma_1|-|\sigma|}{|\sigma_1|}+ \frac{|\sigma_2|-|\sigma|}{|\sigma_2|})=\mathrm{diam}(V_1 \bigcup V_2)(\int\limits_{\sigma_1 \setminus \sigma_2} |\theta_1 -\theta_2| d \mu_H+ \int\limits_{\sigma_2 \setminus \sigma_1} |\theta_1 -\theta_2| d \mu_H) \leq \mathrm{diam}(V_1 \bigcup V_2)\int\limits_{\sigma_1 \bigcup \sigma_2} |\theta_1 -\theta_2| d \mu_H$. This completes the proof.
\end{proof}

Moreover as an easy corollary from statements 3 and 7 we get the proposition

\begin{statement}
The mean point $\mathrm{mp}(V)$ minimizes function which is majorized by barycenter one for any compact $V$ in Hadamard space $H$.
\end{statement}

\section{Behaviour of the Chebyshev center}

\subsection{Reduction to disjoint set of points}

Let us consider properties of the mapping $\mathrm{cheb}: K(H) \to H$.

\begin{statement}
In any finite-dimensional strictly convex metric space $H$ Chebyshev center and radius of the set $V$ can be determined by finite number of points from $V$.
\end{statement}

\begin{proof}
Let us first consider $n=\mathrm{dim} H$ points $x_1, \ldots, x_{n}$ in the intersection $I=S_{r_c(V)} (\mathrm{ch}(V)) \bigcap V$, here $r_c, \mathrm{ch}(V)$ are Chebyshev radius and center of the set $V$. Since $H$ is strictly convex, the set $T=S_{r_c(V)}(x_1) \bigcap \ldots \bigcap S_{r_c(V)}(x_{n})$ is discrete and finite.
Let us denote elements of $T$ by $y_i$, $i=1, \ldots, k$. Then for any $y_i \in T$, $y_i \neq \mathrm{ch}(V)$ there exists a point $z_i \in V$ such that $z_i \not\in B_{r_c(V)}(y_i)$. This holds true since Chebyshev center is unique for any subset of convex metric space. Let us add $z_i$ to the set $I$ and repeat the procedure for any other point of $T$. Thus there exist no more than $2 n -1$ points of $V$ that determine its Chebyshev radius and center.
\end{proof}

\begin{statement}
In any strictly convex infinite-dimensional metric space $H$ both Chebyshev center $\mathrm{cheb}(V)$ (in case such a point exists) and radius $r(V)$ of the set $V$ can be determined by the number of points from $V$ of cardinality less than equal to dimension of $H$.
\end{statement}

\begin{proof}
Consider $V \subset H$.
Since by assumption there exists $\mathrm{cheb}(V)$, we can consider also the ball $B_{\mathrm{cheb}(V)}(r(V))$. Now we must construct $r(V)$-net $V_{r(V)}$ of the set $S_{\mathrm{cheb}(V)}(r(V)) \bigcap V$ and apply statement from \cite{GK}.

Now this net may not be sufficient to determine Chebyshev center of $V$. Then the intersection $I=\bigcap\limits_{x \in V_{r(V)}} B_{x}(r(V))$ contains more than one point. Consider the projection $\pi_i(I)$ of $I$ onto some coordinate line of $H$. Suppose $\pi_i(I)$ consists of more than one point. Then considerations similar to that of the previous statement provide us with the sequence of points $(x^i_k)_{k=1}^{\infty}$ of $V$, such that $\pi_i(I \bigcap B_{x_k^i}(r(V))) \to \pi_i(\mathrm{cheb}(V))$, $k \to \infty$. Since $V$ is compact this sequence possess a converging subsequence $(x^i_{k_l})_{l=1}^{\infty}$, $x^i_{k_l} \to x^i$, $l \to \infty$. Now we simply add this limit point to the set $V_{r(V)}$. The same procedure must be applied to each coordinate line of $H$.

The cardinality of the set of added points is clearly less than equal to the dimension of $H$. This completes the proof.
\end{proof}

\subsection{Behaviour of Chebyshev center}

\begin{statement}
The Chebyshev selection map $\mathrm{cheb}: K(H) \to H$, $V \mapsto \mathrm{cheb}(V)$ is generalised H\"older with power constant $1/2$ \cite{BI, SV} in the Hadamard space $H$.
\end{statement}

\begin{proof}
The proof is done using ideas of \cite{IS, IS2}.
First recall from statement 9 that for any set $V \in K(H)$ its Chebyshev center is uniquely determined by the finite set of points $\sigma(V) \subset \Sigma_{\dim(H)+1}$. 

Let us show now that any variation of some point of $\sigma(V)$ can be represented as combination of the deformations along edges of $co(\sigma(V))$.

The next step is to prove that supremum of the relation $d(\mathrm{cheb}(\sigma), \mathrm{cheb}(\sigma'))/\mathrm{Hd}(\sigma, \sigma')$ is achieved for $\sigma$ and $\sigma'$ being $3$-nets described as follows: $\sigma=\{x, y, z\}$, $\sigma'=\{x, y, z'\}$, here $z \in [x, z']$, $x, y, z \in S_{d(x, y)/2}(m(x, y))$, $x, y, z' \in S_{d(x, z')/2}(m(x, z'))$. 
In order to prove this one must consider first the modulus of convexity of the space of non-positive curvature. It is known \cite{Khamsi} that it is a quadratic function. Thus limit of the relation $d(\mathrm{cheb}(\sigma), \mathrm{cheb}(\sigma'))/\mathrm{Hd}(\sigma, \sigma')$ equals $\infty$ for $z \to y$. 

Let us show now that the relation $d(\mathrm{cheb}(\{x, y, z_1\}), \mathrm{cheb}(\{x, y, z_2\}))/\mathrm{Hd}(\{x, y, z_1\}, \{x, y, z_2\})$ is bounded from above.
Consider $x, y \in H$ and a geodesic ray $\gamma$ starting at $x$. It suffices to show that there exists a constant $L(\gamma, x, y) \in \mathbb{R}^+$ such that for any two points $z, z'$ between $z_1$ and $z_2$. $d(\mathrm{cheb}(\{x, y, z\}), \mathrm{cheb}(\{x, y, z'\})) \leq L d(z, z')$. The only problem here is to determine the behaviour of the function $f(z, z')=d(\mathrm{cheb}(\{x, y, z\}), \mathrm{cheb}(\{x, y, z'\}))/d(z, z')$ for $d(z, z') \to 0$. Assume that for some sequences $z_n$ and $z{'}_n$ $f(z_n, z^{'}_n) \to \infty$. Then there exists a limit point $\gamma(t_0)$ on the geodesic ray $\gamma$. 
Note now that since Chebyshev radius is a Lipschitz function it suffices to estimate $Hd(S_{\mathrm{cheb}(\{x, y, z\})}(r), S_{\mathrm{cheb}(\{x, y, z'\})}(r'))$. Thus $Hd(S_{\mathrm{cheb}(\{x, y, z\})}(r), S_{\mathrm{cheb}(\{x, y, z'\})}(r'))/ d(z, z') \to \infty$. Hence the chord $[y, z]$ of the sphere $S_{\mathrm{cheb}(\{x, y, z\})}(r)$ is tangent to the it. Hence there exists an infinite number of geodesic lines tangent to this chord at point $z$. This is the contradiction with convexity of $H$.


These two facts combined prove that there exist $w, w' \in H$ sufficiently close to $y$ for any $(z_1, z_2) \in [z, z'] \times [z, z']$ such that both of the following claims hold true

1) $w \in [x, w']$, $x, y, w \in S_{d(x, y)/2}(m(x, y))$, $x, y, w' \in S_{d(x, w')/2}(m(x, w'))$;

2) $d(\mathrm{cheb}(\{x, y, z_1\}), \mathrm{cheb}(\{x, y, z_2\}))/\mathrm{Hd}(\{x, y, z_1\}, \{x, y, z_2\}) \leq$\\ $\leq d(\mathrm{cheb}(\{x, y, w\}), \mathrm{cheb}(\{x, y, w'\}))/\mathrm{Hd}(\{x, y, w\}, \{x, y, w'\})$.

The last estimate then can be derived from consideration of \cite{Khamsi} on convexity modulus of non-positively curved spaces.

The general case of $n$-nets can be analysed similarly. 
\end{proof}

\begin{example}
The behaviour of the Chebyshev center for the spaces of Busemann nonpositive curvature is not as easily described as in the case of Hadamard manifolds and surely does not possess the same nice properties. 

1. $H=(\mathbb{R}^2, d =\|\cdot\|_p)$, $p>2$. Then for the point $(1, 0)$ of the unit sphere $S_0(1)$ we get the convexity modulus equal to $\varepsilon^p$. Thus the mapping $\mathrm{cheb}: (\Sigma_3(\mathbb{R}^2), \mathrm{Hd}) \to (\mathbb{R}^2, \|\cdot\|_p)$ is H\"older in the neighbourhood of $\sigma \in \Sigma_2$, $\sigma=\{(-1, 0), (1, 0)\}$ with the power coefficient equal to $1/p$.

2. $H=(\bigotimes\limits_{i=1}^{\infty}\mathbb{R}^2, \|\cdot\|)$, here $\|x\|=(\sum\limits_{i=1}^{\infty}(|x_{2i}|^{i+2}+|x_{2i -1}|^{i+2})^{\frac{2}{i+2}})^2$. Note first that this space is uniform convex in any direction, so Chebyshev center is unique \cite{Gar}. Then the convexity of this space in any direction is not bounded from below by $\varepsilon^p$ for any $p \in \mathbb{N}$. Thus Chebyshev selection map is not even H\"older one.
\end{example}

The only thing we can be sure at is that the behaviour of the Chebyshev center is the best in $CAT(0)$ of all convex metric spaces.

\begin{statement}
Let $H$ be strictly convex space of Busemann nonpositive curvature.
Then the upper bound for the power of the generalised H\"older map $\mathrm{cheb}: B(H) \to H$ is $1/2$.
\end{statement}

\begin{proof}
Consider the construction involving $\sigma=\{x, y, z\}$ and $\sigma'=\{x, y, z'\}$ from the previous statement.

This is the direct consequence of the convexity properties of spheres from $H$.
 
The only obstacle is estimation of the distance not between points on the spheres but the centers of them. Assume that $d(m(x, z'), m(x, y)) = o(d(z, z'))$. Then there exist two geodesic lines passing through $x$ in the same direction. This contradicts convexity of $H$. Hence for any $y \in S_{1}(x)$ there exists $L_x>0$ such that $d(m(x, z'), m(x, y))/  d(z, z') \geq L_x$. 
\end{proof}

Note that uniform convexsity of the space is crucial for H\"older behaviour of the Chebyshev center.

\begin{example}
Let us describe behaviour of Chebyshev center for $(\mathbb{R}^n, \|\cdot \|_1)$. 

First note that to prove that the mapping $\mathrm{cheb}$ is Lipschitz it suffices to find upper estimate for the distance between angle points of the hyperplane subset $\mathrm{cheb}\{x, y\}$ for arbitrary pair $x, y \in \mathbb{R}^n$. In order to do this consider a pair of points $x=(x_1, \ldots, x_n)$ and $y=(y_1, \ldots, y_n)$. Let us denote by $t=(t_1, \ldots, t_n)$ points of the set $\mathrm{cheb}(\{x, y\})$. Then the angle points of $\mathrm{cheb}(\{x, y\})$ have coordinates $t_1, \ldots, t_i, \ldots, t_n$, here $t_j$ equals either $t_{j, m}=\min \{x_j, y_j\}$ or $t_{j, M}=\max \{x_j, y_j\}$ and $|t_i-x_i|=\frac{\sum\limits_{j=1}^{n}|x_j -y_j|-2\sum\limits_{j\neq i}|t_j -t_{j, m}|}{2}$. Consider $I=\{j\in\{1, \ldots, n\}| t_j =t_{j, m}\}$. Assume without loss of generality that $x_i \leq y_i$. Now let $x'$ and $y'$ of $\mathbb{R}^n$ be such that both $\|x-x'\|_1$ and $\|y-y'\|_1$ are less or equal to $\varepsilon>0$.
Then 
$$
\|t-t'\|_1=\sum\limits_{j=1}^{n}|t_j - t_j^{'}|=\sum\limits_{j\neq i}|t_j-t_j^{'}|+|t_i -t_i^{'}|
\leq
2\varepsilon+ 
$$
$$
+
\frac{|\sum\limits_{j=1}^{n}|x_j -y_j|-2\sum\limits_{j\neq i}|t_j -t_{j, m}|-\sum\limits_{j=1}^{n}|x^{'}_j -y^{'}_j|-2\sum\limits_{j\neq i}|t^{'}_j -t^{'}_{j, m}|+2(x_i -x_i^{'})|}{2}
\leq
$$
$$
\leq
2 \varepsilon + \frac{|\sum\limits_{j\not\in I}|x_j -y_j| -\sum\limits_{j \in I}|x_j -y_j|-\sum\limits_{j\not\in I}|x^{'}_j -y^{'}_j|+ \sum\limits_{j\in I}|x^{'}_j -y^{'}_j|+2(x_i -x_i^{'})|}{2}
\leq
4\varepsilon.
$$
Hence since Chebyshev center of the subset of $(\mathbb{R}^n, \|\cdot \|_1)$ is determined by at most $n$ pairs of points, $L \leq 4n$.

Next let us show that there exists also a lower estimate for Lipschitz constant of $\mathrm{cheb}$. In order to show this consider Chebyshev center of the sets 
$$
\sigma_1=\{(0, \ldots, 0), (1, \ldots, 1), (0, 1, \ldots, 1), (1, 0, \ldots, 0), \ldots, (1, \ldots, 1, 0, 1), (0, \ldots, 0, 1, 0)\}
$$ 
and 
$$
\sigma_2=\{(\varepsilon/n, \ldots, \varepsilon/n), (1+\varepsilon/n, \ldots, 1+\varepsilon/n), (-\varepsilon/n, 1+\varepsilon/n, \ldots, 1+\varepsilon/n), 
$$
$$
(1-\varepsilon/n, \varepsilon/n, \ldots, \varepsilon/n),
\ldots, 
$$
$$
(1+\varepsilon/n, \ldots, 1+\varepsilon/n, -\varepsilon/n , 1+\varepsilon/n), (\varepsilon/n, \ldots, \varepsilon/n, 1-\varepsilon/n, \varepsilon/n)\}.
$$ 
Then $d(\mathrm{cheb}(\sigma_1), \mathrm{cheb}(\sigma_2)) \geq (n-1) \varepsilon$. Thus $L \geq n-1$.
\end{example}

The second similar example ---$l^{\infty}$-space --- was presented in \cite{IS2}.

\section{Lipschitz selection for $\Sigma_n(H)$.}

In this section we construct Lipschitz selections for the convex hulls of no more than $n$ points of metric space $H$. 

\subsection{Construction based on the mapping $\mathrm{mp}$}
Let us now construct the mapping $c$ from the set of all $n$-nets $\Sigma_n \subset \mathbb{R}^n$ to $\mathbb{R}^n$, such that $c$ is Lipschitz with constant $4$. This means that $d(c(\sigma), c(\sigma'))\leq \alpha(\sigma, \sigma')$, here $\alpha$ is Hausdorff metric.

Let us first consider this problem locally.

The construction is inductive.

It is clear that $c(\sigma)=x$ for any $\sigma =\{x\} \in \Sigma_1$. This map obviously is Lipschitz.

To make our construction correct we must assume also that 
\begin{equation}
\forall \sigma, \sigma' \in \Sigma_n, \sigma \setminus\{x_n\}=\sigma'\setminus\{x'_n\},\;
d(x_n, x'_n) \leq \varepsilon \Rightarrow d(c(\sigma), c(\sigma')) \leq \varepsilon/2.
\end{equation}

This obviously holds true for $\Sigma_2$ and $c(\{x_1, x_2\})=m(x_1, x_2)$.

Assume that the construction is valid for any $\Sigma_l$, here $1\leq l\leq n-1$.

Consider $\sigma \in \Sigma_n$. 

We put $c(\sigma)=b(\sigma)$ if for any point $x_i \in \sigma$, $d(x_i, co(\sigma\setminus\{x_i\}))\geq \mathrm{diam}(\sigma)/2$.

Assume now that there exist points $x_1, \ldots, x_k$, $k\leq n$, such that the distances $d(x_i, co(\sigma\setminus\{x_i\}))$ are less than $\mathrm{diam}(\sigma)/2$.

Let us construct the set $C$ consisting of points $c_i=c(\sigma \setminus x_i)$, $i=\overline{1, k}$ and
consider the barycenter $b(C)$.  

Consider now the segment $[b(\sigma), b(C)]$. 
We define $c(\sigma) \in [b(\sigma), b(C)]$ to be the point which divides $[b(\sigma), b(C)]$ in relation $d(b(C), c) : d(b(\sigma), b(C)) = \mathrm{min}\{d(x_i, co(\sigma\setminus \{x_i\}))\}_{i=\overline{1,k}}/\mathrm{diam}(\sigma) : 1$.

Now we must verify condition (1) for the set $\Sigma_n$.

By assumption points $c_i=c(\sigma \setminus \{x_i\})$ depend on $x_i$ in Lipschitz way with the constant $L_{n-1}$.

The mapping $c$ is continuous by construction. 

The only thing we must show is that its Lipschitz constant exists. Let us estimate this constant as one of the combination of the mappings, namely, shifts of the points of $C$ and path along the segment itself. 
Then the constant $L_n=L_{n-1}$(shift of the point $m(C)$) $+1$ (path along the segment; it does not depend on the number of shifted point but only on the minimal shift) $ + 1$ (shift of the point $m(\sigma)$). Thus $L_n=2+L_{n-1}$, hence the constant exists.

Note now that Hausdorff metric locally coincides with Fedorchuk one \cite{Fed}. Recall also that the latter metric is inner, thus the estimate $L_n$ for Lipschitz constant of $c: C(X) \to X$ is globally true for the space of $n$-nets endowed with Fedorchuk metric.

\subsection{Selection based on the mapping $\mathrm{cheb}$}

Now assume as in the previous section that there exists a Lipshitz selection $l: \Sigma_{n-1}(H) \to H$ with Lipschitz constant $L_{n-1}$. Let us expand this mapping to $\Sigma_n(H)$. Consider Chebyshev center $\mathrm{cheb}(\sigma)$ of the set $\sigma \in \Sigma_n(H)$. It is known that the selection $\mathrm{cheb}:\Sigma_n(H) \to H$ is generalised H\"older one with H\"older constant $H_n \mathrm{diam}^{1/2}(\sigma)$ and power coefficient $1/2$. Thus as in the first section we consider the set $\{\sigma_1, \ldots, \sigma_n\}$ of $n-1$ subnets of $\sigma$. Then there exists a set $\sigma_l=\{l(\sigma_1), \ldots, l(\sigma_n)\}$. Consider barycenter $m(\sigma_l)$ and connect it with Chebyshev center of $\sigma$ by segment. We define $l(\sigma) \in [m(\sigma_l), \mathrm{cheb}(\sigma)]$ as the point which divides $[m(\sigma_l), \mathrm{cheb}(\sigma)]$ in relation $d(l(\sigma), m(\sigma_l)) : d(\mathrm{cheb}(\sigma), l(\sigma)) = \mathrm{min}\{d^{1/2}(x_i, co(\sigma\setminus \{x_i\}))\}_{i=\overline{1,k}}/(\mathrm{diam}(\sigma) \max\{1, H_n\}): 1 - \mathrm{min}\{d^{1/2}(x_i, co(\sigma\setminus \{x_i\}))\}_{i=\overline{1,k}}/\mathrm{diam}(\sigma) \max\{1, H_n\})$. This construction provides us with the desired selection $l: \Sigma_n(H) \to H$ with Lipschitz constant equal to $1+L_{n-1}+3/2 L_{n-1}$.

It seems that the last construction can be easily extended to the space $K(H)$. That is, one may take into consideration $n$-nets defining Chebyshev center of the set $V \subset K(H)$ together with all their possible $n-1$-subnets. Nevertheless, there exists a natural obstruction, namely, a decent measure on the so-called perimeter or boundary of the set $V$. The author still does not know the solution of this problem.

\subsection{Construction of the Lipschitz selection for the subsets of convex hulls of finite sets}

The most natural generalisation of the construction given in previous paragraph is as follows:


\end{document}